\documentclass{scrartcl}
\usepackage[T1]{fontenc}
\usepackage[utf8]{inputenc}
\usepackage{amsmath, amsthm, amssymb}
\usepackage{dsfont}
\usepackage{graphicx}
\usepackage[format=plain]{caption}%ben?tigt, um mit dem Befehl \caption*{} eine Bildunterschrift zu setzen, die nicht mit Abbildung X eingeleitet wird. 'normal'=\setcapindent{0}
\usepackage{capt-of}%ben?tigt, um 'figure' und 'table' nebeneinander zu setzen, siehe 'http://www.faqs.org/faqs/de-tex-faq/part6/' Punkt 6.1.7, aufgerufen am 01.08.11.
\usepackage{subcaption}
\usepackage{bm}
\renewcommand{\d}{\,\mathrm{d}} %Differenzial in Integralen
\renewcommand{\phi}{\varphi}
\newcommand{\ul}{\underline}
\newcommand{\ol}{\overline}

\numberwithin{equation}{section}
\theoremstyle{definition} 
\theoremstyle{definition}

\title{Simulation and optimal control of the Williams-Otto process using Pyomo}

\author{Jochen Schmid$^{1}$, Katrin Teichert$^{1}$, \\Moncef Chioua$^{2}$, Thorsten Schindler$^{2}$, Michael Bortz$^{1}$\\  
\small $^1$ Fraunhofer Institute for Industrial Mathematics (ITWM), 67663 Kaiserslautern, Germany \\
\small $^2$ ABB AG, Corporate Research Germany, %Wallstadter Straße 59,
68526 Ladenburg, Germany}  

\date{}

\begin{document}

\maketitle

\begin{abstract}
\small{ \noindent 
We illustrate the advantages the high-level open-source software package Pyomo has 
%can provide 
%for the chemical engineer 
in rapidly setting up and solving dynamic simulation and optimization problems. In order to do so, we use the example of the Williams-Otto process. We show how to simulate the process dynamics using the collocation method and the IPOPT solver provided by Pyomo. We also discuss waste minimization and yield maximization as two examplary process optimization problems. And finally, we present and compare two approaches to setpoint tracking: one based on proportional-integral feedback control and one based on optimal open-loop control.
%
%In recent years, a number of sophisticated dynamic optimization software tools such as Pyomo, ACADO and CasADi have emerged, which opens up new possibilities for chemical engineers. With these frameworks, it becomes increasingly easy to set up and solve ODE and DAE models of chemical processes and then to augment these models in order to solve optimal control or data fitting problems. %by introducing an objective function and free controls/parameters. 
%In this paper, we illustrate the advantages of Pyomo on
%the prototypical example of the Williams-Otto process. We show how to simulate the process dynamics %by solving the underlying ODE system 
%using the collocation method and the IPOPT solver provided by Pyomo. We also discuss waste minimization and yield maximization as two examplary process optimization problems. And finally, we present and compare two approaches to setpoint tracking: one based on proportional-integral feedback control and one based on optimal open-loop control. 
}
\end{abstract}

{ \small \noindent 
%\emph{AMS Subject Classification (2010):} 
%\\
}

\section{Introduction}

Virtualization permeates many aspects of chemical engineering nowadays. Computer models allow process optimization without the costs associated with experiments. And conversely, parameter estimation and machine learning methods allow experimental data to be fed back into computer models, increasing their reliability. 
\smallskip

Chemical process models can be static (steady-state) or dynamic. Capturing the dynamic aspect becomes increasingly important if the process is supposed to adapt to external changes, such as fluctuations in supply, demand or energy cost. A prominent example of a method based on a dynamic model of the process is model-predictive control~\cite{Al}, \cite{FiAlBi07}, \cite{GrPa}. In case its underlying dynamic model is accurate enough, model-predictive control has the potential for smoother and thus more efficient adaptations compared to simple feedback loops like proportional-integral-derivative controllers. %Also, model-predictive control can handle complex control scenarios like  multivariable systems with constraints on both control and output variables, systems with significant time delays, or sytems with inverse dynamic response. With proportional-integral-derivative control, by contrast, handling complex control scenarios usually requires ad hoc (costly development-phase) solutions. 
In spite of these drawbacks, proportional-integral-derivative controllers are still dominant in the process industries due to their simplicity~\cite{ViVi12} -- but in recent years, with the advent of faster and faster optimization solvers (Table~2 and~3 of~\cite{AnDi+19}), linear and nonlinear model-predictive control has found more and more industrial applications as well~\cite{FiAlBi07}. In particular, model-predictive control will probably be used more systematically in online optimization~\cite{GrKrRa01}, which could replace the standard industrial real-time optimization solutions based on static process models that are currently offered by many vendors and used by the process industry.  
%In recent years, however, with the advent of faster and faster optimization solvers (Table~2 and~3 of~\cite{CasADi}), model-predictive control has found more and more industrial applications as well~\cite{} and it will probably be used more systematically in online optimization~\cite{GrKrRa00}.  
\smallskip
 
A dynamic model can be the basis for a process optimization (optimal control) or data fitting task (parameter estimation). In either case, we obtain a so-called dynamic optimization problem~\cite{Biegler10}
\begin{gather}
\min J(\bm{x},\bm{u},\bm{p}) \label{eq:opt-contr-general}\\
\bm{f}(\dot{\bm{x}}(t),\bm{x}(t),\bm{u}(t),\bm{p},t) = 0 \qquad (t \in [0,t_f]) \label{eq:DAE-constraint} \\
\bm{g}(\bm{x}(t),\bm{u}(t),\bm{p},t) \le 0 \qquad (t \in [0,t_f]). \label{eq:path-constraint-general-DAE}
\end{gather}
where we distinguish between state trajectories $\bm{x}$, control profiles $\bm{u}$ to be optimized, and time-independent finite-dimensional parameters $\bm{p}$ to be optimized. Sometimes, the final time $t_f$ is an optimization variable as well, in which case it has to be added to the goal function $J$'s list of arguments, but we will not consider that in the present paper. In chemical engineering applications, the control profiles $\bm{u}$ can be the temperature profile in a reactor or the reflux ratio profile in a batch distillation column~\cite{BhBi96}, \cite{CeBi98}. And the time-independent parameters $\bm{p}$ will typically be model parameters -- like reaction rates~\cite{TjBi91} or equlibrium model parameters~\cite{BoHe+19} -- or design parameters -- like reactor length, diameter, inlet temperature and composition in reactor design, for instance. 
In the simplest case, the inequality constraints~\eqref{eq:path-constraint-general-DAE} just take the form of box constraints
\begin{align} \label{eq:box-constraints-general}
\ul{\bm{x}} \le \bm{x}(t) \le \ol{\bm{x}}, \qquad \ul{\bm{u}} \le \bm{u}(t) \le \ol{\bm{u}}, \qquad \ul{\bm{p}} \le \bm{p} \le \ol{\bm{p}}  \qquad (t \in [0,t_f])
\end{align}
with appropriate lower and upper bounds for the state and control variables $\bm{x}, \bm{u}$ and the parameters $\bm{p}$. (As usual, for vectorial quantities, the above inequalities are to be understood componentwise.) Clearly, inequality constraints of the form~\eqref{eq:path-constraint-general-DAE} also comprise equality constraints -- in the case of box constraints~\eqref{eq:box-constraints-general}, for instance, one just has to choose identical lower and upper bounds. 
In many cases, the differential-algebraic constraints~\eqref{eq:DAE-constraint} can be separated into an ordinary differential equation for the states and additional algebraic equations, so that the optimal control problem~\eqref{eq:opt-contr-general}-\eqref{eq:path-constraint-general-DAE} then takes the simpler form 
\begin{gather}
\min J(\bm{x},\bm{u},\bm{p}) \\
\dot{\bm{x}}(t) = \bm{f}(\bm{x}(t),\bm{u}(t),\bm{p},t) \qquad (t \in [0,t_f]) \label{eq:ODE-constraint} \\
\bm{g}(\bm{x}(t),\bm{u}(t),\bm{p},t) \le 0 \qquad (t \in [0,t_f]).  \label{eq:path-constraint-general-ODE}
\end{gather}

Approaches to solve dynamic optimization problems as above can be distinguished according to how strongly the procedure of solving the differential equation is intertwined with the procedure of improving the goal function value. There is a whole spectrum of different approaches  ranging from methods with weakly intertwined to methods with strongly  intertwined solution and optimization procedures. %See, for instance, \cite{Biegler10} (Section~8.6) for a compact overview.  %weak/strong intertwining
Also, in recent years, a number of sophisticated software tools to solve dynamic optimization problems have emerged like Pyomo~\cite{Pyomo}, \cite{Ni+18}, ACADO~\cite{HoFeDi11} and CasADi~\cite{AnDi+19}, for instance. 
See, for instance, \cite{Biegler10} (Section~8.6) for a compact overview over the spectrum of solution approaches for dynamic optimization problems.  
\smallskip 

At the weakly intertwined end of the spectrum lies the so-called sequential approach. %At the end of the spectrum where solution and optimization are least intertwined
In this approach, a single step in the controls $\bm{u}$ and parameters $\bm{p}$ towards improvement of the objective function (for example in the direction of steepest descent) alternates with a complete solution of the differential equation, for example using a Runge-Kutta method.
\smallskip
 
At the strongly intertwined end of the spectrum lies the so-called simultaneous approach and, in particular, the collocation method used in Pyomo~\cite{Pyomo}, \cite{Ni+18}. In this approach, both states and controls are discretized simultaneously. The differential equation is then transformed into a system of equations by use of orthogonal collocation. The algebraic constraints are also transformed to a finite equation system. After this transformation, the discretized states and their derivatives, the discretized controls and the finite-dimensional parameters are all treated as optimization parameters within one large  nonlinear optimization problem.
\smallskip
 
Somewhere in the middle between these two ends of the spectrum lies the multiple shooting~\cite{BoPl84}, \cite{DiBo+06} approach used in ACADO~\cite{HoFeDi11}. In that approach, the time horizon $[0,t_f]$ of the states is partitioned. Similarly to the sequential approach, a single step towards improvement of the objective alternates with solving the differential equations for all partition intervals. However, the so-called matching conditions, which enforce continuity of the states from one partition interval to the next, are part of the optimization problem. This is achieved by treating the initial values for each partition interval as additional optimization parameters. Therefore, similarly to the collocation approach, the solution of the differential equation over the whole time horizon is obtained alongside the improvement of the objective, and only at convergence the solution $(\bm{x}^*,\bm{u}^*,\bm{p}^* )$ is guaranteed to globally solve the differential equation.

\section{Pyomo} \label{sect:pyomo}

Pyomo~\cite{Ha+11}, \cite{Pyomo}, \cite{Ni+18} was originally developed by researchers in the Center for 
Computing Research at Sandia National Laboratories and is a COIN-OR 
project. The software can be used to solve dynamic optimization problems as formulated in~\eqref{eq:opt-contr-general}-\eqref{eq:DAE-constraint}. In particular, with the pyomo.dae module it is possible to solve optimal control problems with differential-algebraic or even partial differential-algebraic equations (DAE or PDAE) as constraints~\cite{Ni+18}. 
\smallskip

Pyomo is based on Python, a popular high-level programming language. Aside from a small number of Pyomo-specific commands, no specific modelling language has to be learnt to start using it. The user defines the dynamic optimization problem in its original form~\eqref{eq:opt-contr-general}-\eqref{eq:DAE-constraint} using constructs provided by Pyomo. For example, the time domain $[0,t_f]$ is introduced as a continuous set, and the states and their derivatives are identified as such by specific commands. Pyomo features automatic differentiation, so the user does not have to specify functions to return the derivatives of the objective function $J$ nor of the constraints $\bm{f}$ and $\bm{g}$. 
\smallskip

%Additionally, all discretizations of the differential equations and controls are done automatically in Pyomo. 
Additionally, Pyomo provides two automatic discretization transformations: orthogonal collocation and finite differences. See  \cite{Ni+18}, \cite{Pyomo} or the Pyomo website for algorithmic details. Both transformations take a user-defined continuous model and discretize it automatically, that is, transform the dynamic optimization problem into a standard non-linear program composed of algebraic functions and equations. In particular, no manual discretziations are necessary, which removes many potential sources of errors.
\smallskip

Finally, the discretized optimization model is handed to the IPOPT solver included in Pyomo. IPOPT is a powerful state-of-the-art interior-point solver~\cite{WaBi06}, \cite{WaBi05a}, \cite{WaBi05b}. It should be noticed that the discretization transformation leads to a largely sparse algebraic system. Pyomo ensures that IPOPT makes use of this property when solving the discretized problem.

\section{Simulation of the Williams-Otto process}

\subsection{Williams-Otto process}

%After having introduced our simulation tool of choice, it is now time to introduce the process to which we want to apply it, namely the Williams-Otto process. 
In this section, we introduce the process on which we illustrate the simulation framework from Section~\ref{sect:pyomo}, namely the Williams-Otto process. 
The Williams-Otto process is an exemplary chemical process which comprises several %the most 
common process units from chemical engineering and is therefore extensively studied in the literature~\cite{WiOt60}, \cite{LoVaHoDi12}, \cite{HaMa12}, \cite{CaKw16}, \cite{Biegler97}. %has therefore been extensively studied in the literature to the present day 

\begin{figure}[htbp]
\centering
\includegraphics[width=10cm]{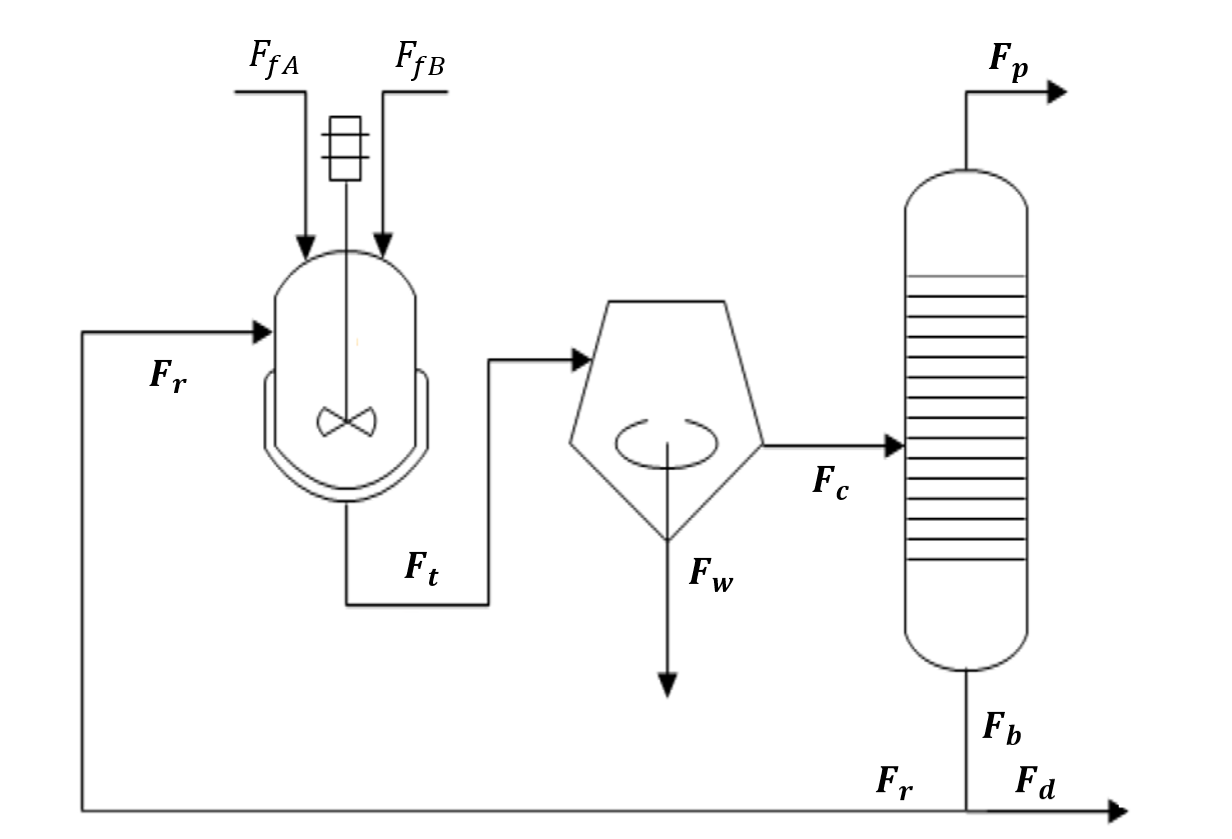}
\caption{Williams-Otto process flowsheet (adapted from~\cite{CaKw16})}
\label{fig:WO-flowsheet}
\end{figure}

It consists of four subsystems: a continuously stirred tank reactor, a decanter (centrifuge), a distillation column (separator) and a splitter. See Figure~\ref{fig:WO-flowsheet} for a schematic flowsheet of the process.  
Two substances $A$ and $B$ are fed into the reactor (with mass streams $F_{fA}$ and $F_{fB}$ respectively), where they interact according to the following reaction system:
\begin{gather}
A + B \longrightarrow C \label{eq:reaction-1}\\
C + B \longrightarrow P + E \label{eq:reaction-2}\\
P + C \longrightarrow G \label{eq:reaction-3}
\end{gather} 
%The reaction constants depend on the temperature T inside the reactor.
A compound mass stream $\bm{F}_t =(F_{tA},F_{tB},F_{tC},F_{tP},F_{tE},F_{tG})$ of total amount $\mu := F_{tA} + \dotsb + F_{tG}$ then leaves the reactor. In the decanter, the waste product $G$ is removed from the process through a waste stream %$\bm{F}_w$:
\begin{align} \label{eq:Fw_G}
F_{wG} = F_{tG}.
\end{align}
In the distillation column, part of the desired product $P$ is distilled with a fixed efficiency, namely 
\begin{align} \label{eq:Fp_P}
F_{pP} = F_{cP}-0.1 F_{cE} = F_{tP}-0.1 F_{tE}
\end{align}
where $\bm{F}_c$ is the input stream to the distillation column. The bottom stream $\bm{F}_b$ of the distillation column is partly withdrawn at a ratio $\eta$ (stream $\bm{F}_d$) and partly recycled back to the reactor at the ratio $1-\eta$ (stream $\bm{F}_r$). %The bottom stream $\bm{F}_b$ of the distillation column is split: namely, it is partly withdrawn at a ratio $\eta$ (stream $\bm{F}_d$) and partly recycled back to the reactor at the ratio $1-\eta$ (stream $\bm{F}_r$). 
In total, there are five controls: the feed streams $F_{fA}$ and $F_{fB}$, the reactor temperature $T$, the total mass stream $\mu$ leaving the reactor, and the split fraction $\eta$. Also, there are six substances $A, B, C, P, E, G$ that interact in the reactor, and their masses inside the reactor are denoted by $m_A, \dots, m_G$.  
\smallskip

Applying the rules of reaction kinetics and using component mass balances around each of the four subsystems, one can derive the following system of ordinary differential equations for the masses $m_A, \dots, m_G$ of the species in the reactor:  
\begin{align}
\dot{m}_A &= F_{f A} + ((1-\eta) \mu - \mu) \frac{m_A}{m}  - k_1(T) \frac{m_A m_B}{V} \label{eq:component-mass-balance-in-terms-of-m-A}\\
\dot{m}_B &=  F_{f B} + ((1-\eta) \mu - \mu) \frac{m_B}{m} - k_1(T) \frac{m_A m_B}{V} - k_2(T) \frac{m_B m_C}{V} \\
\dot{m}_C &=  ((1-\eta) \mu - \mu) \frac{m_C}{m}  + 2 k_1(T) \frac{m_A m_B}{V} - 2 k_2(T) \frac{m_B m_C}{V}  - k_3(T) \frac{m_C m_P}{V} \\
\dot{m}_E &=  ((1-\eta) \mu - \mu) \frac{m_E}{m}  +  2 k_2(T) \frac{m_B m_C}{V}\\
\dot{m}_P &=  0.1 (1-\eta) \mu \frac{m_E}{m} - \mu \frac{m_P}{m}  + k_2(T) \frac{m_B m_C}{V} - 0.5 k_3(T) \frac{m_C m_P}{V} \\
\dot{m}_G &=  - \mu \frac{m_G}{m}  + 1.5 k_3(T) \frac{m_C m_P}{V}, \label{eq:component-mass-balance-in-terms-of-m-G}
\end{align}
In these equations, $m := m_A + \dotsb + m_G$ and $V := m_A/\rho_A^{0} + \dotsb + m_G/\rho_G^{0}$ denote the total mass and the total volume of the mixture in the reactor, respectively, with equal pure-substance densities $\rho_A^{0}, \dots, \rho_G^{0} = 50 \, \mathrm{lb/ft^3}$ (Appendix~I of~\cite{WiOt60}). Also, the reaction rates $k_i$ are given by 
\begin{align} \label{eq:k_i'}
k_i(T) = \frac{a_i}{\rho} \, \mathrm{exp}(-b_i/T) \qquad (i \in \{1,2,3\})
\end{align}
where $a_i$, $b_i$, $\rho$ are the constants from~\cite{Biegler97} (Section~8.2.1):
\begin{gather} \label{eq:a_i,b_i,rho-Biegler}
a_1 := 5.9755 \cdot 10^9 \, \mathrm{h}^{-1}, \quad a_2 := 2.5962 \cdot 10^{12} \, \mathrm{h}^{-1}, \quad a_3 := 9.6283 \cdot 10^{15} \, \mathrm{h}^{-1},\\
\rho := 50 \, \mathrm{lb}/\mathrm{ft}^3, \quad b_1 := 12000 \, \mathrm{^\circ R}, \quad b_2 := 15000 \, \mathrm{^\circ R}, \quad b_3 := 20000 \, \mathrm{^\circ R} 
\end{gather}
It is through these reaction rates that the reactor temperature $T$ -- which is one of our controls -- influences the process. 
In the entire paper, we adopt the common unit conventions from the Williams-Otto literature~\cite{WiOt60}, \cite{CaKw16}, \cite{Biegler97}, that is, masses are measured in $\mathrm{klb}$, time is measured in $\mathrm{h}$, mass flows are measured in $\mathrm{klb/h}$, %volumes are measured in $1000 \, \mathrm{ft}^3$
and temperatures are measured in $^\circ \mathrm{R}$ (degrees R\'{e}aumur).
\smallskip

We will often refer to~\eqref{eq:component-mass-balance-in-terms-of-m-A}-\eqref{eq:component-mass-balance-in-terms-of-m-G} as the Williams-Otto equations in what follows. In view of the bilinear mass terms in the numerators and the total mass and total volume terms in the denominators on the right-hand sides of~\eqref{eq:component-mass-balance-in-terms-of-m-A}-\eqref{eq:component-mass-balance-in-terms-of-m-G}, the Williams-Otto equations are clearly nonlinear. 
With the abbreviations
\begin{gather} 
\bm{x}(t) := (m_A(t), m_B(t), m_C(t), m_E(t), m_P(t), m_G(t)), \\
\bm{u}(t) := (F_{f A}(t), F_{f B}(t), T(t), \mu(t), \eta(t))
\end{gather}
for the state and control variables, the Williams-Otto equations assume the standard form 
\begin{align}
\dot{\bm{x}}(t) = \bm{f}(\bm{x}(t),\bm{u}(t))
\end{align}
of a nonlinear time-invariant finite-dimensional control system. It can be shown that for every given initial state $\bm{x_0}$ and every given (bounded) control profile $\bm{u}$, the corresponding initial value problem
\begin{align}
\dot{\bm{x}}(t) = \bm{f}(\bm{x}(t),\bm{u}(t)) \qquad \text{and} \qquad \bm{x}(0) = \bm{x_0}
\end{align}
has a unique (local) solution in the sense of Carath\'{e}odory. See~\cite{Sontag} (Lemma~2.6.2) or~\cite{BressanPiccoli} (Section~3.2), for instance. 
All streams of the Williams-Otto process -- except for the feed streams, of course -- can be directly calculated from the solution of~\eqref{eq:component-mass-balance-in-terms-of-m-A}-\eqref{eq:component-mass-balance-in-terms-of-m-G} and the tank outlet stream $\mu$: for example, 
\begin{align} \label{eq:Ft-in-terms-of-m-and-mu}
F_{t A} = \mu \frac{m_A}{m}, \quad \dots \quad F_{t G} = \mu \frac{m_G}{m}.
%F_{t S} = \mu \frac{m_S}{m} \qquad (S \in \{A, \dots, G\}).
\end{align}

\subsection{Simulation examples} \label{sect:simul-examples}

As has been pointed out in Section~\ref{sect:pyomo}, implementing the Williams-Otto equations in Pyomo is fairly easy and convenient. %See Section~\ref{sect:pyomo} and~\cite{Ni+18}, \cite{Pyomo}. 
In this subsection, we present two very simple simulation examples for the Williams-Otto process performed in Pyomo. In both examples, we start from the initial state $\bm{x_0} := (10, 1, 0,0,0,0)$ and we impose a piecewise constant control input $\bm{u}$ with just one jump (step change) at the time $t^* := 100$, that is, 
\begin{align*}
\bm{u}(t) = \bm{u}^* \quad (t \le t^*) \qquad \text{and} \qquad \bm{u}(t) = \bm{u}^{**} \quad (t > t^*).
\end{align*}
Also, the initial part $\bm{u}^*$ of the control is the same in both examples, namely
\begin{align} \label{eq:u*}
F_{fA}^* := 10, \quad F_{fB}^* := 20, \quad T^* := 580, \quad \mu^* := 129.5, \quad \eta := 0.2. 
\end{align} 

\begin{figure}[htbp]
\centering
\includegraphics[width=\textwidth]{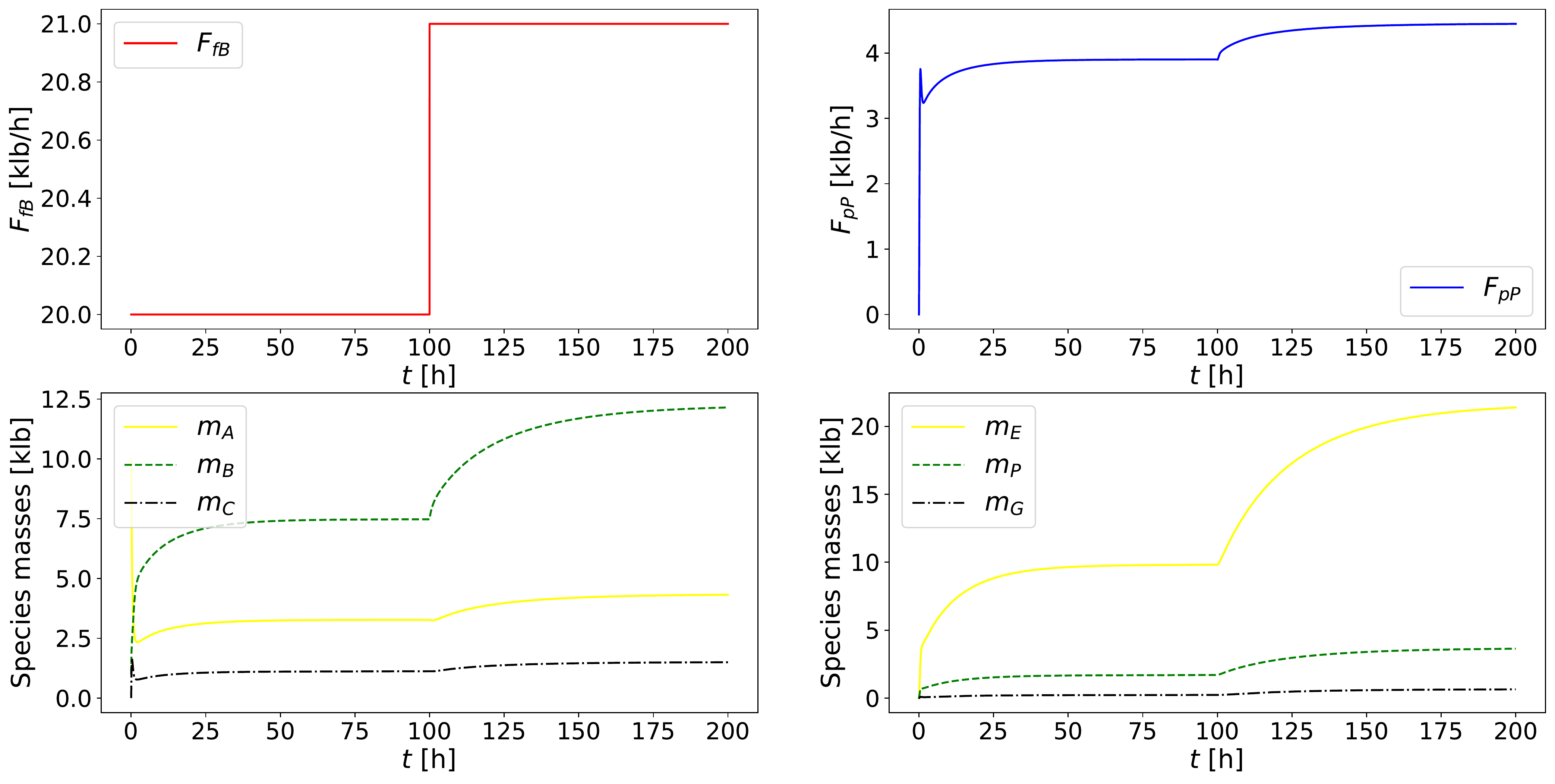}
\caption{Step change in $F_{fB}$, response of $F_{pP}$ (top row), and resulting state trajectories (bottom row)}
\label{fig:simul-Ff_B-Fp_P}
\end{figure}

\begin{figure}[htbp]
\centering
\includegraphics[width=\textwidth]{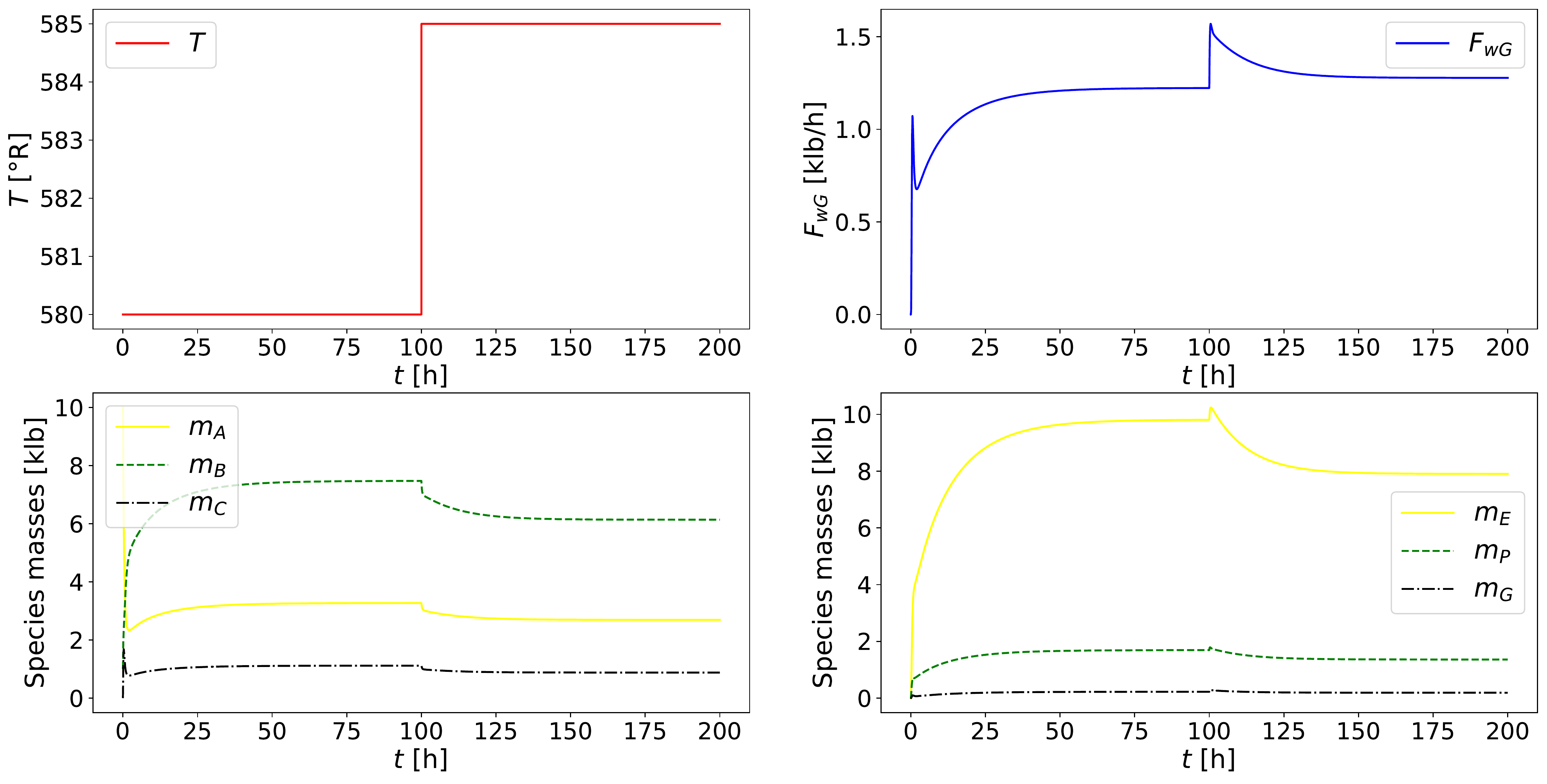}
\caption{Step change in $T$, response of $F_{wG}$ (top row), and resulting state trajectories (bottom row)}
\label{fig:simul-T-Fw_G}
\end{figure}

\noindent In the first example (Figure~\ref{fig:simul-Ff_B-Fp_P}), the step change is in the control variable $F_{fB}$, namely from $F_{fB}^* = 20$ to $F_{fB}^{**} := 21$. And in the second example (Figure~\ref{fig:simul-T-Fw_G}), the step change is in the $T$-profile, from $T^* = 580$ to $T^{**} := 585$. It takes just about $1$ CPU second and $7$ IPOPT iterations to perform the above simulation examples (chosen number of finite elements and number of collocation points per finite element: $1500$ and $3$, respectively). 
We see that before and long after the jump, the state trajectories become flat (Figure~\ref{fig:simul-Ff_B-Fp_P} and~\ref{fig:simul-T-Fw_G} (bottom row)). In other words, at those times, the Williams-Otto process is in steady-state operation. %is in a steady state.

\section{Setpoint tracking by means of proportional-integral control}

An important task in process engineering is setpoint tracking, that is, to steer a process from a steady state to another steady state in such a way that a chosen process quantity of interest -- the so-called setpoint quantity --  asymptotically approaches a user-specified setpoint. In our Williams-Otto example, typical setpoint quantities are the product stream $F_{pP}$ from the head of the distillation column and the waste stream $F_{wG}$ from the bottom of the decanter. 
\smallskip

A standard approach to perform setpoint tracking is to use proportional-integral controllers. Such a controller -- just like any other controller -- is based on an appropriately designed feedback relation between the setpoint quantity $y$ in question and an appropriate corresponding control quantity $u$. In formulas, the feedback relation of a proportional-integral controller looks as follows:
\begin{align} \label{eq:pi-control,general}
u(t) = u^{bias} + K^{prop}(y^{sp}-y(t)) + K^{int} \int_{t^*}^t (y^{sp}-y(\tau))\d \tau  
\end{align}  
for $t > t^*$, that is, $t^*$ is the time where the feedback loop between setpoint and control quantity is activated. %where $t^*$ is the time at which the controller is switched on. 
We immediately see from this formula that the control action $u(t)$ is the larger, the larger the current and the cumulated deviation 
\begin{align} \label{eq:prop-and-int-term}
y^{sp}-y(t), \qquad \int_{t^*}^t (y^{sp}-y(\tau)) \d \tau
\end{align}
of our setpoint quantity $y$ from the user-specified setpoint $y^{sp}$ are. (We tacitly used here the fact that the parameters $K^{prop}$ and $K^{int}$ have the same sign, see Section~\ref{sect:controller-parameters} below.) In other words, the farther our setpoint quantity is %we are 
away from the desired setpoint, the more actively our controller works to push the setpoint quantity towards its setpoint. %the more active is our controller in order to counteract and diminish the deviation.  
We also see from the above formula~\eqref{eq:pi-control,general} that as soon as the setpoint $y^{sp}$ has been reached sufficiently well, the control action $u$ does not change appreciably anymore. %remains constant just because the proportional and the integral term in~\eqref{eq:prop-and-int-term} are then  equal to zero or equal to a constant, respectively. 
While in this case the proportional term in~\eqref{eq:pi-control,general}  is equal to zero, the integral term in~\eqref{eq:pi-control,general} is equal to a generally nonzero constant. It is this nonzeroness that enables a proportional-integral controller to steer a process from a steady state $(\bm{x}^*,\bm{u}^*)$ to a different steady state $(\bm{x}^{**},\bm{u}^{**})$ with $\bm{u}^{**} \ne \bm{u}^*$.  
\smallskip

Commonly, the parameters $u^{bias}$, $K^{prop}$, $K^{int}$ of a proportional-integral controller are called the control bias, the proportional and the integral controller gain, respectively. Also, any pair $(u,y)$ of a control and a setpoint quantity related by the proportional-integral feedback relation~\eqref{eq:pi-control,general} is called a control channel. 
Although proportional-integral controllers are particularly well-suited for linear processes, they can also be used to perform small setpoint changes (around a fixed operation point) in nonlinear processes like the Williams-Otto process. %they can also be used for setpoint tracking of nonlinear processes like the Williams-Otto process. 
In fact, we use proportional-integral controllers for the control channels 
\begin{align} \label{eq:contr-channels-WO}
(\mu, m), \qquad (F_{fA},F_{tP}), \qquad (F_{fB},F_{pP}), \qquad (T,F_{wG})
\end{align}
of the Williams-Otto process. Similar control channels are proposed %used 
in~\cite{CaKw16}.
It is relatively easy to implement the proportional-integral controllers corresponding to the control channels~\eqref{eq:contr-channels-WO} in Pyomo. All one has to do is to complement the Williams-Otto equations~\eqref{eq:component-mass-balance-in-terms-of-m-A}-\eqref{eq:component-mass-balance-in-terms-of-m-G} %by the feedback relation~\eqref{eq:pi-control,general} with $(u,y)$ being one of the control channels from~\eqref{eq:contr-channels-WO}. 
by the respective feedback relations
\begin{align} 
\mu(t) &= \mu^{bias} + K^{prop}_1 (m^{sp}-m(t)) + K^{int}_1 \int_{t^*_1}^t (m^{sp}-m(\tau))\d \tau  \label{eq:pi-control-1,WO}\\
F_{fA}(t) &= F_{fA}^{bias} + K^{prop}_2 (F_{tP}^{sp}-F_{tP}(t)) + K^{int}_2 \int_{t^*_2}^t (F_{tP}^{sp}-F_{tP}(\tau))\d \tau  \label{eq:pi-control-2,WO}\\
F_{fB}(t) &= F_{fB}^{bias} + K^{prop}_3 (F_{pP}^{sp}-F_{pP}(t)) + K^{int}_3 \int_{t^*_3}^t (F_{pP}^{sp}-F_{pP}(\tau))\d \tau  \label{eq:pi-control-3,WO}\\
T(t) &= T^{bias} + K^{prop}_4 (F_{wG}^{sp}-F_{wG}(t)) + K^{int}_4 \int_{t^*_4}^t (F_{wG}^{sp}-F_{wG}(\tau))\d \tau \label{eq:pi-control-4,WO}
\end{align}  
for $t$ larger than the respective channel's activation time $t_i^*$. In view of~\eqref{eq:Fw_G}, \eqref{eq:Fp_P} and \eqref{eq:Ft-in-terms-of-m-and-mu}, the streams $F_{tP}$, $F_{pP}$ and $F_{wG}$ explicitly depend on the tank outlet stream $\mu$. And therefore the controllers~\eqref{eq:pi-control-2,WO}-\eqref{eq:pi-control-4,WO}, when used in conjunction with~\eqref{eq:pi-control-1,WO}, implicitly also contain the feedback relation~\eqref{eq:pi-control-1,WO}.

\subsection{Controller tuning} \label{sect:controller-parameters}

In order to obtain a controller with a good tracking performance, it is essential to carefully choose the controller parameters $u^{bias}, K^{prop}, K^{int}$~\cite{O'Dwyer12}. In this subsection, we explain how this parameter tuning is done. %how this can be done. %how we did that. 
As a first step, one has to choose a design level of operation, that is, a typical steady state $(\bm{x}^*,\bm{u}^*)$ around which one wishes to operate the considered process. As our design level of operation we choose the steady state $(\bm{x}^*,\bm{u}^*)$ from Section~\ref{sect:simul-examples}, that is, 
\begin{align} \label{eq:design-level-of-operation}
\bm{x}^* := (3.27, 7.47,1.12,9.81, 1.69, 0.22)
\end{align}
and $\bm{u}^*$ is given by~\eqref{eq:u*}. 
Accordingly, we then choose the bias terms as
\begin{align*}
\mu^{bias} := \mu^* = 129.5, \quad F_{fA}^{bias} := F_{fA}^{*} = 10, \quad F_{fB}^{bias} := F_{fB}^{*} = 20, \quad T^{bias} := T^* = 580
\end{align*}  
to guarantee bumpless control. 
%According to this design level of operation, we choose the bias terms as
%\begin{align*}
%\mu^{bias} := \mu^* = 129.5, \quad F_{fA}^{bias} := F_{fA}^{*} = 10, \quad F_{fB}^{bias} := F_{fB}^{*} = 20, \quad T^{bias} := T^* = 580
%\end{align*}  
%(thereby guaranteeing bumpless control). 
In order to determine the controller gains $K_i^{prop}$, $K_i^{int}$ we perform step tests for our four control channels. A step test means that, starting from our design level of operation, we abruptly step up or down the control quantity from $u^*$ to a new value $u^{**}$ (just like we did in Section~\ref{sect:simul-examples}). We then observe how our setpoint quantity responds to the step in $u$ and wait until the process  settles to a final steady state $(\bm{x}^{**},\bm{u}^{**})$ (or at least until the setpoint quantity of interest does not change anymore). 
\smallskip

In the case of the control channels $(\mu, m), (F_{fA},F_{tP}),  (F_{fB},F_{pP})$, the step response is of first order with deadtime $0$ (Figure~\ref{fig:simul-Ff_B-Fp_P}). Accordingly, we can use the Skogestad rules of internal model control~\cite{Sk06}, \cite{SkGr12} 
\begin{align}
K^{prop} := \frac{1}{K^{proc}} \qquad \text{and} \qquad K^{int} := \frac{t^{proc}}{K^{proc}}
\end{align}
(in the version of moderate tuning). In these relations, $K^{proc}$ is the process static gain of the considered control channel, which relates the total eventual change $\Delta y$ in the setpoint quantity to the percentage jump $\Delta_{rel} u$ in the control quantity. %which expresses how much the setpoint quantity changes relative to the percentage jump in the control quantity, 
Also, $t^{proc}$ is the process time constant of the considered control channel, which is defined as the timespan it takes for the setpoint quantity to perform $63 {\%}$ of its total change $\Delta y$. In short,  
\begin{align*}
K^{proc} := \Delta y /\Delta_{rel} u \qquad \text{and} \qquad 
y(t^*+t^{proc}) = y^* + 0.63 \cdot \Delta y,
\end{align*}
where $\Delta y := y^{**}-y^* = y(t_f) - y(t^*)$ and $\Delta_{rel} u := (u^{**}-u^*)/(\ol{u}-\ul{u}) = (u(t_f)-u(t^*))/(\ol{u}-\ul{u})$ with $\ul{u}, \ol{u}$ being typical lower and upper bounds for $u$. 
We thus find 
\begin{gather*}
K_1^{prop} = -0.002, \quad K_1^{int} = -40.659, \quad K_2^{prop} = 0.053, \quad K_2^{int} = 1.504, \\
K_3^{prop} = 0.069, \quad K_3^{int} = 1.282
\end{gather*}
In the case of the control channel $(T,F_{wG})$, the step response is not of first order (Figure~\ref{fig:simul-T-Fw_G}). Accordingly, we cannot apply the Skogestad rules but instead have to tune the controller gains by hand. We found the following values to yield a good tracking performance: 
\begin{align*}
K_4^{prop} = 0.143, \quad K_4^{int} = 0.1
\end{align*}

\subsection{Simulation examples} \label{sect:pi-results}

\begin{figure}[htbp]
\centering
\includegraphics[width=\textwidth]{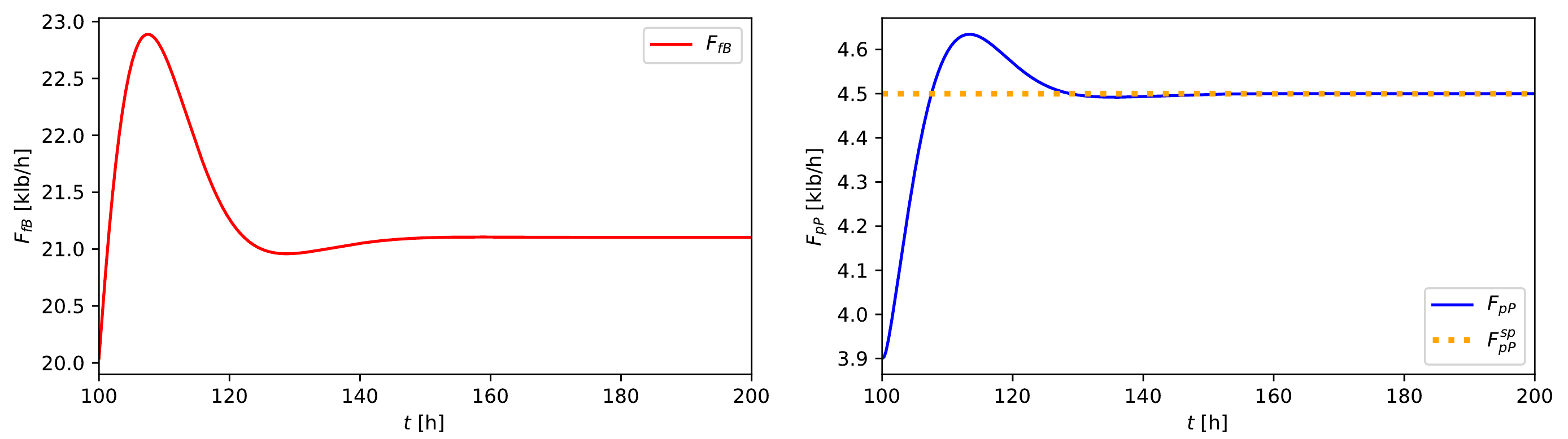}
\caption{Setpoint tracking for $F_{pP}$ by proportional-integral control}
\label{fig:setpt-Fp_P-pi}
\end{figure}

\begin{figure}[htbp]
\centering
\includegraphics[width=\textwidth]{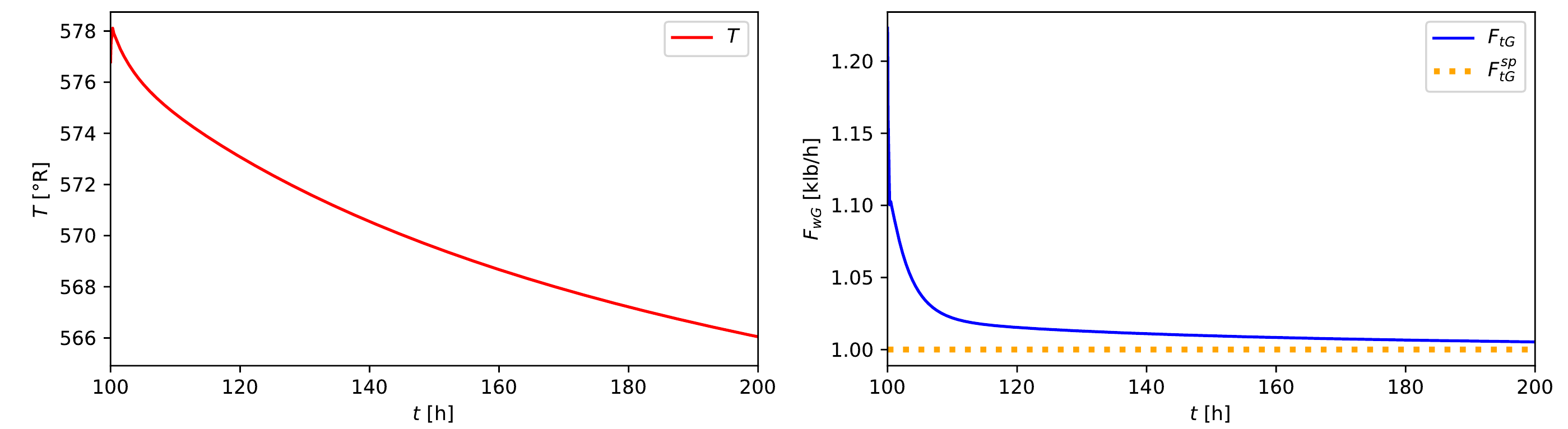}
\caption{Setpoint tracking for $F_{wG}$ by proportional-integral control}
\label{fig:setpt-Fw_G-pi}
\end{figure}

In this subsection, we present two setpoint tracking examples for the Williams-Otto process: one with the control channel $(F_{fB},F_{pP})$ and one with the control channel $(T,F_{wG})$. In both examples, we initialize  the process as in the examples from Section~\ref{sect:simul-examples}. And then at time $t^* = 100$, we switch on our proportional-integral controller~\eqref{eq:pi-control-3,WO} or~\eqref{eq:pi-control-4,WO} with some exemplary setpoint value $F_{pP}^{sp}$ or $F_{wG}$, respectively. See Figure~\ref{fig:setpt-Fp_P-pi}  and \ref{fig:setpt-Fw_G-pi} with the setpoints depicted in orange. 
Since proportional-integral controllers are linear, they generally exhibit good tracking performance only for small setpoint changes around the chosen design level of operation. 
%As expected, %the setpoints are reached after a reasonble time
%We found that the setpoint tracking works well as long as the setpoints do not lie too far away from the values the setpoint quantity has before the controller is switched on. 
%\smallskip
%
We can also use several of our controllers~\eqref{eq:pi-control-1,WO}-\eqref{eq:pi-control-4,WO} simultaneously. When doing so, however, we are not as free anymore in the choice of setpoints as we were in the case of a single control channel: for example, not every pair $(F_{pP}^{sp},F_{wG}^{sp})$ of setpoints for $F_{pP}$ and $F_{wG}$ can be reached by using~\eqref{eq:pi-control-3,WO} and~\eqref{eq:pi-control-4,WO} simultaneously.

\section{Sample optimal control problems}

Apart from setpoint tracking, another important task in process engineering is optimal control, that is, to optimize certain quality measures of a process by choosing optimally shaped control profiles. In our Williams-Otto example, typical quality measures are the total amount of waste and the total amount of product accumulated over the whole process duration. In general, the optimal control problems considered here read as follows: %an optimal control problem reads as follows: %minimize or maximize a certain goal function $J(x,u)$ under constraints of the form
\begin{align} \label{eq:opt-contr-goal}
\text{minimize or maximize } J(\bm{x},\bm{u}) 
\end{align}
($J$ being the desired process quality measure) subject to constraints of the form
\begin{gather}
\dot{\bm{x}}(t) = \bm{f}(\bm{x}(t),\bm{u}(t)) \qquad (t \in [0,t_f]) \qquad \text{and} \qquad \bm{x}(0) = \bm{x_0} \label{eq:ode-constraint} \\
\bm{g}(\bm{x}(t),\bm{u}(t)) \le 0 \qquad (t \in [0,t_f]) \label{eq:opt-contr-general-ineq-constraint}
\end{gather}
(differential-equation constraints and inequality constraints). It is easy to implement such optimal control problems in Pyomo. All one has to do is to complement the process equations~\eqref{eq:ode-constraint} by the following specifications: one has to
\begin{itemize}
\item specify the goal function $J$ one wishes to minimize or maximize
\item specify which control variables should be the free optimization variables and which shape the corresponding control profiles should have (general spline, piecewise constant, or constant)
\item specify inequality constraints of the form~\eqref{eq:opt-contr-general-ineq-constraint} on the control and state variables, like box constraints on the controls or path constraints on state or output trajectories, for instance.
\end{itemize}

\subsection{Waste minimization} \label{sect:waste-minimization}

\begin{figure}[htbp]
\centering
\includegraphics[width=\textwidth]{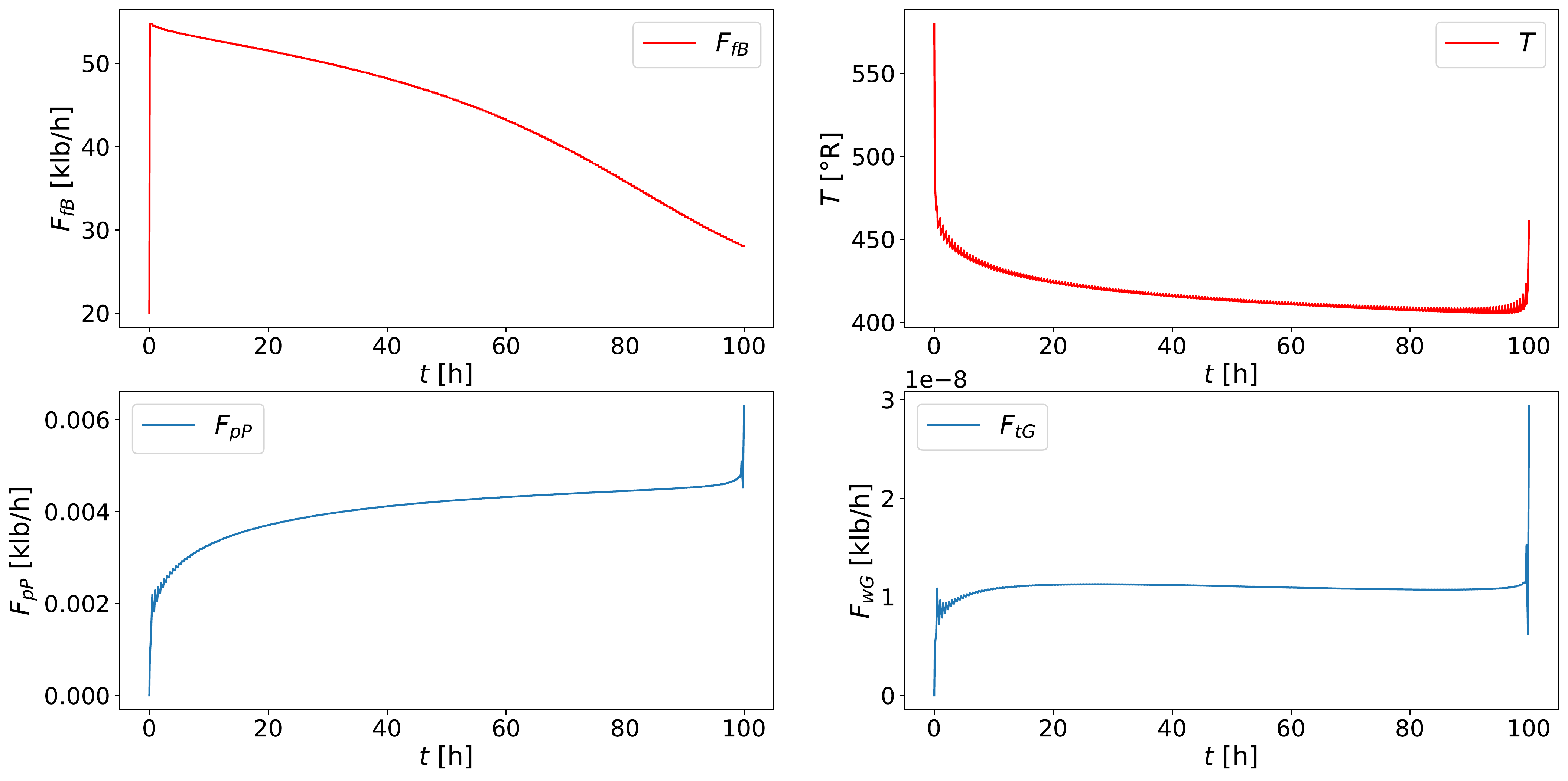}
\caption{Control profiles minimizing the total waste (top row), and resulting stream profiles (bottom row)}
\label{fig:waste-minimize}
\end{figure}

In this subsection, we briefly discuss the waste minimization problem, that is, the goal function $J$ to be minimized is given by
\begin{align}
J(\bm{x},\bm{u}) := \int_0^{t_f} F_{wG}(\tau) \d\tau,
\end{align}
where the waste stream $F_{wG}$ depends on $(\bm{x},\bm{u})$ according to~\eqref{eq:Fw_G} and~\eqref{eq:Ft-in-terms-of-m-and-mu} and where as total process time we choose $t_f = 100$. 
As our free optimization variables we choose the feed stream $F_{fB}$ (profile type: piecewise constant) and the reactor temperature $T$ (profile type: general spline), as initial state $\bm{x_0}$ we choose the steady state~\eqref{eq:design-level-of-operation}, and as constraints we choose the following box constraints %box constraints for these control variables
\begin{align*}
F_{fA} = 10, \quad 0 \le F_{fB} \le 56, \quad 200 \le T \le 800, \quad \mu = 129.5, \quad \eta = 0.2
\end{align*}
After $257$ iterations and $16.7$ CPU seconds, the IPOPT solver then finds the control profiles from Figure~\ref{fig:waste-minimize} (top row) to be optimal (chosen number of finite elements and collocation points per finite element: $200$ and $3$, respectively). %and the corresponding minimal value of our total waste goal function is  
With these control profiles, the value of our total waste goal function is basically $0$ but, on the other hand, the total yield of our desired product is very small as well (Figure~\ref{fig:waste-minimize} (bottom row)).  %the total amount of desired product accumulated at the head of the distillation column 
We therefore consider an appropriate yield maximization problem next. 

\subsection{Yield maximization} \label{sect:yield-maximization}

\begin{figure}[htbp]
\centering
\includegraphics[width=\textwidth]{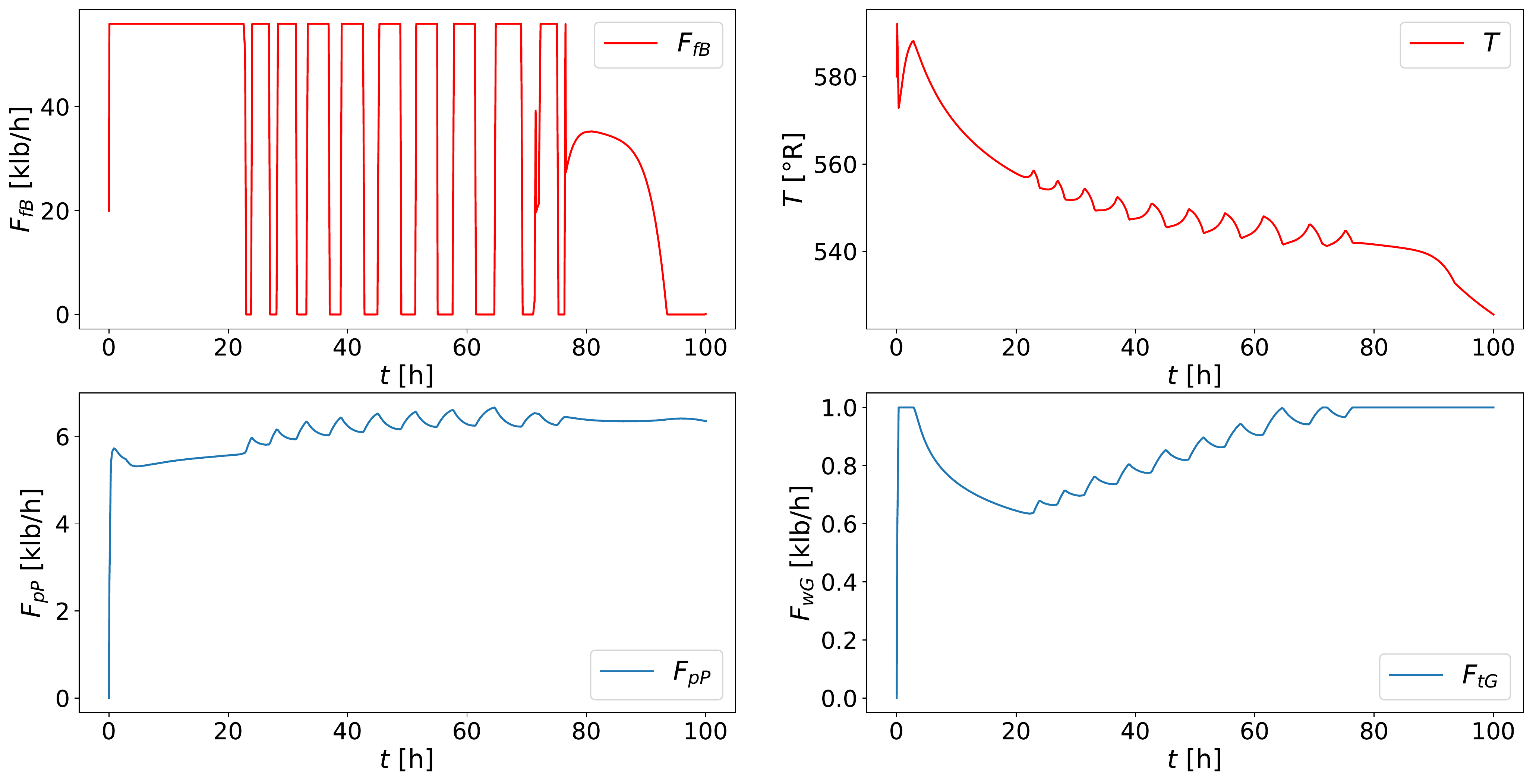}
\caption{Control profiles maximizing the total yield under the waste constraint~\eqref{eq:waste-constraint} (top row), and resulting stream profiles (bottom row)}
\label{fig:yield-maximize-path-constr}
\end{figure}

In this subsection, we discuss the yield maximization problem under a waste constraint, that is, the goal function $J$ to be maximized is given by
\begin{align} \label{eq:yield-goal-fct}
J(\bm{x},\bm{u}) := \int_0^{t_f} F_{pP}(\tau) \d\tau,
\end{align}
where the product stream $F_{pP}$ depends on $(\bm{x},\bm{u})$ according to~\eqref{eq:Fp_P} and~\eqref{eq:Ft-in-terms-of-m-and-mu} and where as total process time we choose $t_f = 100$. 
As our free optimization variables we choose the feed stream $F_{fB}$ (profile type: general spline) and the reactor temperature $T$ (profile type: general spline), as initial state $\bm{x_0}$ we choose the steady state~\eqref{eq:design-level-of-operation}, and as constraints we choose the following box constraints %for our control variables 
\begin{align*}
F_{fA} = 10, \quad 0 \le F_{fB} \le 56, \quad 200 \le T \le 800, \quad \mu = 129.5, \quad \eta = 0.2
\end{align*}
and, additionally, the following path constraint on the waste stream:
\begin{align} \label{eq:waste-constraint}
F_{wG}(t) \le 1 \qquad (t \in [0,t_f]).
\end{align}
After $1048$ iterations and $75.2$ CPU seconds, the IPOPT solver then finds the control profiles from Figure~\ref{fig:yield-maximize-path-constr} (top row) to be optimal (chosen number of finite elements and collocation points per finite element: $200$ and $3$, respectively). %and the corresponding minimal value of our total waste goal function is  
With these control profiles, our yield goal function takes the value $611.3$. 
We see that the optimal $F_{fB}$-profile is bang-bang~\cite{BressanPiccoli}, \cite{Sontag} for a long period of time, that is, switching between the lower and upper bound for $F_{fB}$. Such a bang-bang behavior is not very popular in industrial practice as it can lead to an accelerated wear of the actuators. It is possible to reduce the bang-bang behavior of $F_{fB}$, without substantially decreasing the achieved yield value, 
\begin{itemize}
\item by adding an appropriate upper-bound constraint on the integrated $F_{fB}$-control action $\int_0^{t_f} F_{fB}(\tau) \d \tau$ (hard constraint)
\item or by adding the integrated $F_{fB}$-control action $\int_0^{t_f} F_{fB}(\tau) \d \tau$ with an appropriate weight to the original goal function~\eqref{eq:yield-goal-fct} (soft constraint).
\end{itemize}

\noindent We do not go into detail here because we achieved an even better bang-bang suppression (and, at the same time, even better total yield and total waste values) by the combined optimization approach presented next.

%DETAILLIERTERE AUSFUEHRUNGEN ZUR BANG-BANG REDUKTION
%\smallskip
%
%We see that the optimal $F_{fB}$-profile is bang-bang~\cite{BressanPiccoli}, \cite{Sontag} for a long period of time, that is, switching between the lower and upper bound for $F_{fB}$. Such a bang-bang behavior is not very popular in industrial practice as it can lead to an accelerated wear of the actuators. We can reduce the bang-bang behavior of $F_{fB}$ by  complementing the above optimal control problem by an upper bound constraint on the integrated $F_{fB}$-control action like, for instance,
%\begin{align}
%\int_0^{t_f} F_{fB}(\tau) \d \tau \le 2700.
%\end{align}
%With this additional constraint, the bang-bang behavior is substantially reduced while the optimal value of our yield goal function decreases not very much (namely to $577.7$). %deteriorates not very much (namely from $611.3$ to $577.7$). 
%Alternatively, we can reduce the bang-bang behavior of $F_{fB}$ by adding the integrated $F_{fB}$-control action (with an appropriate weight) to the original goal function~\eqref{eq:yield-goal-fct} as a soft constraint. With the new goal function
%\begin{align}
%J(x,u) := \int_0^{t_f} F_{pP}(\tau) \d\tau - 0.05 \int_0^{t_f} F_{fB}(\tau) \d \tau
%\end{align}
%for instance, the bang-bang behavior is substantially reduced while the value of the yield decreases only a little (namely to $599.7$). 

\subsection{Combined waste and yield optimization}
%\subsection{Comined waste minimization and yield maximization}

\begin{figure}[htbp]
\centering
\includegraphics[width=\textwidth]{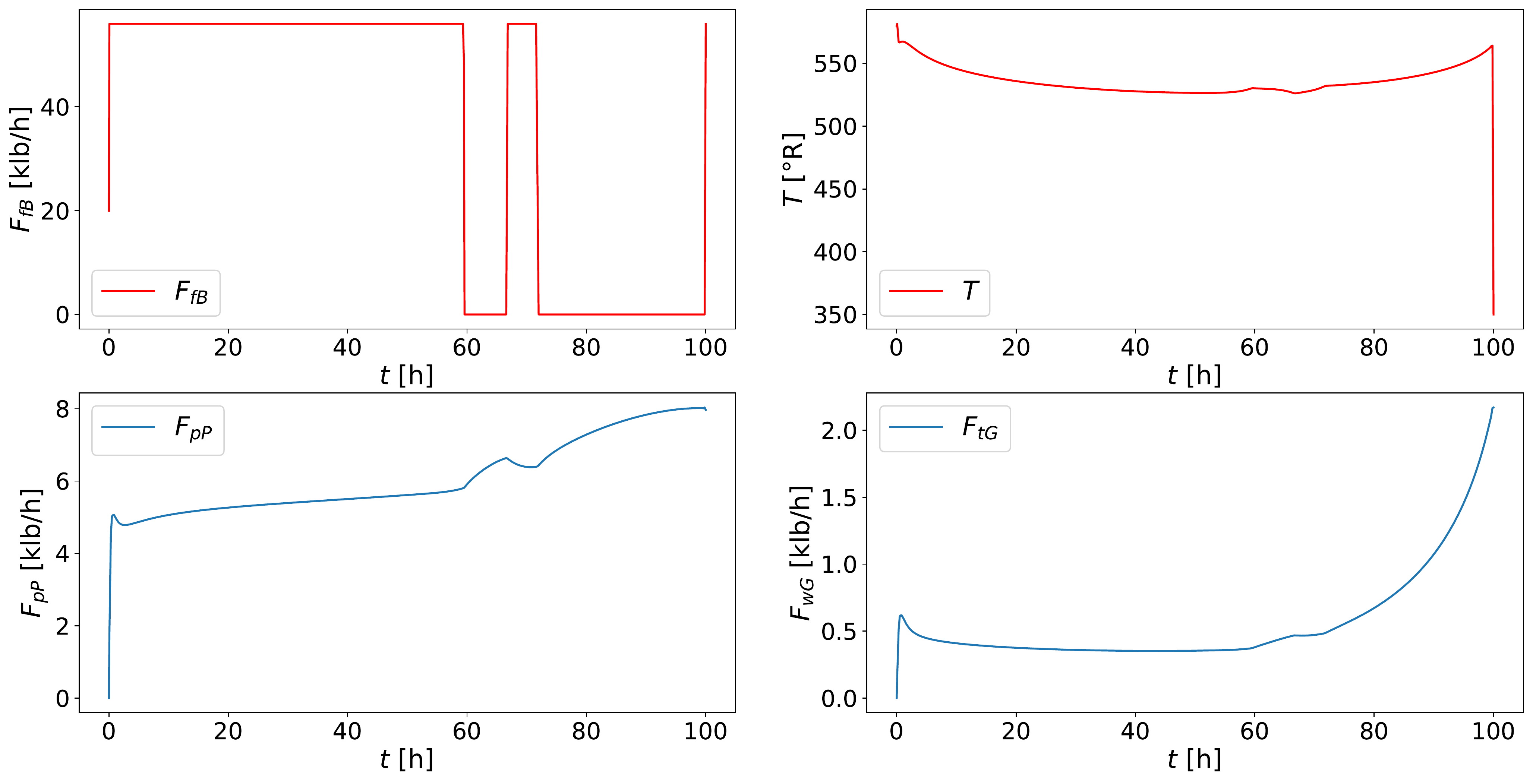}
\caption{Control profiles maximizing the weighted sum~\protect\eqref{eq:weighted-sum-yield-vs-waste} of total yield and total waste (top row), and resulting stream profiles (bottom row)}
\label{fig:combined-waste-and-yield-optimization}
\end{figure}

While in the previous two subsections we optimized waste and yield separately, we now want to optimize these conflicting quality measures  simultaneously. In order to do so, we maximize the goal function $J$ given by
\begin{align} \label{eq:weighted-sum-yield-vs-waste}
J(\bm{x},\bm{u}) := \alpha \int_0^{t_f} F_{pP}(\tau) \d\tau - \beta \int_0^{t_f} F_{wG}(\tau) \d \tau
\end{align}  
with appropriately chosen weights $\alpha, \beta >0$. We choose the same final time $t_f$, the same initial state $\bm{x_0}$, and the same free optimization variables and constraints as in Section~\ref{sect:waste-minimization}. Choosing equal weights $\alpha = 1$ and $\beta=1$, we then  find the control profiles from Figure~\ref{fig:combined-waste-and-yield-optimization} (top row) to be optimal after $388$ IPOPT iterations and $22.9$ CPU seconds (chosen number of finite elements and collocation points per finite element: $200$ and $3$, respectively). %With equal weights $\alpha = 1$ and $\beta=1$, the IPOPT solver then finds the control profiles from Figure~\ref{fig:yield-vs-waste} (top row) to be optimal after $388$ iterations and $22.9$ CPU seconds (chosen number of finite elements and collocation points per finite element: $200$ and $3$, respectively).
With these optimal control profiles, the yield is almost as large as in Section~\ref{sect:yield-maximization} (namely $608.3$ compared to $611.3$) while the total amount of waste is significantly smaller than in Section~\ref{sect:yield-maximization} (namely $56.5$ compared to $86.3$). A nice extra feature is that the strongly oscillatory bang-bang behavior in the optimal $F_{fB}$-profile from Section~\ref{sect:yield-maximization} is  drastically reduced here.

\subsection{Setpoint tracking by means of optimal control}

In the previous section, we discussed the standard approach to setpoint tracking based on feedback control, namely proportional-integral control. In this subsection, we present an alternative approach which is based on optimal control. %Since the setpoint tracking problem is about steering a process to a steady-state setpoint, it is natural to formulate it ...
In that approach, we formulate the setpoint tracking problem -- which, recall, is about steering a process quantity $y$ to a steady-state setpoint $y^{sp}$ -- as an optimal control problem~\eqref{eq:opt-contr-goal}-\eqref{eq:opt-contr-general-ineq-constraint} where the goal function $J$ is given by  
\begin{align} \label{eq:goal-fct-setpoint-tracking}
J(\bm{x},\bm{u}) := \int_{t^*}^{t_f} (y(\tau)-y^{sp})^2 \d \tau + \alpha \dot{\bm{x}}(t_f)^2
\end{align}
and to be minimized. In this expression, $t^*$ is the time instant at which the user wants to start the setpoint tracking and $\alpha > 0$ is a weight parameter. (As usual, the process should be in a steady state at this time $t^*$.)
%and where this goal function is to be minimized. 
Clearly, by minimizing the first part of the goal function~\eqref{eq:goal-fct-setpoint-tracking} we ensure that our setpoint quantity $y$ approaches its setpoint $y^{sp}$ as well and as quickly as possible, while by minimizing the second part of~\eqref{eq:goal-fct-setpoint-tracking} we ensure that the process eventually is as close to a steady state as possible.  
\smallskip

\begin{figure}[htbp]
\centering
\includegraphics[width=\textwidth]{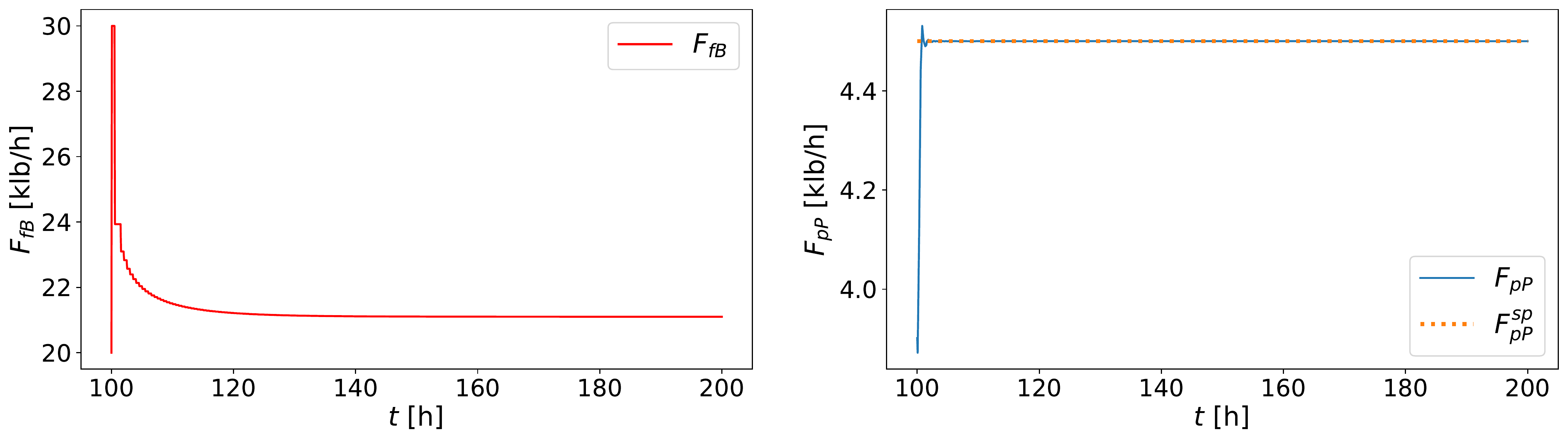}
\caption{Setpoint tracking for $F_{pP}$ by optimal control}
\label{fig:setpt-Fp_P-opt-contr}
\end{figure}

\begin{figure}[htbp]
\centering
\includegraphics[width=\textwidth]{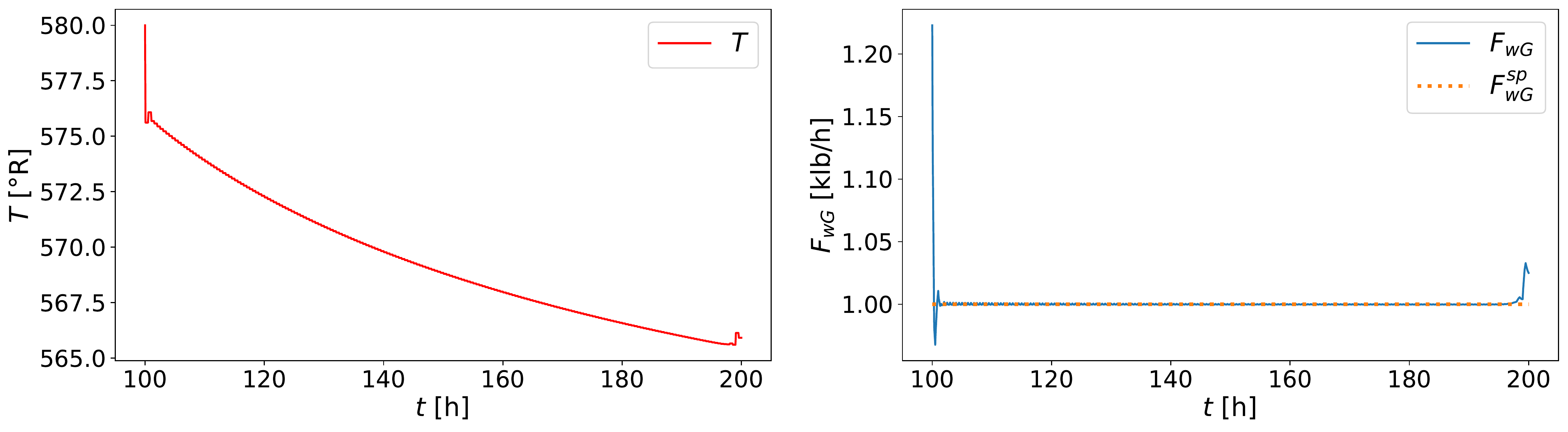}
\caption{Setpoint tracking for $F_{wG}$ by optimal control}
\label{fig:setpt-Fw_G-opt-contr}
\end{figure}

We apply this optimization-based approach of setpoint tracking to the same two examples as in the feedback-based approach before (Section~\ref{sect:pi-results}), that is, as setpoint quantities we choose $F_{pP}$ and $F_{wG}$, respectively, and as corresponding optimal control variables we choose, respectively, $F_{fB}$ and $T$. Also, we choose the same time $t^*$ and the same initial state $\bm{x_0}$ and initial control profile $\bm{u}|_{[0,t^*]} \equiv \bm{u}^*$ %before the setpoint is set. 
as in Section~\ref{sect:pi-results}. And finally, we choose $t_f = 200$ and $\alpha = 1$.
Comparing the results (Figures \ref{fig:setpt-Fp_P-opt-contr} and~\ref{fig:setpt-Fw_G-opt-contr}) with our results from the feedback-based approach (Figures \ref{fig:setpt-Fp_P-pi} and \ref{fig:setpt-Fw_G-pi}), we see that with the optimization-based approach the setpoints are reached much more quickly, namely within less than $10 \%$ ($F_{pP}$ case) or $3 \%$ ($F_{wG}$ case) of the time needed in the feedback-based approach. %and that the driving control quantities oscillate less heavily (smaller amplitudes). 
\smallskip

%We finally point out that 
Additionally, the above optimization-based approach works just as well %for the simultaneous tracking of two setpoint quantities 
for tracking two setpoint quantities simultaneously, like $F_{pP}$ and $F_{wG}$ (instead of just a single one). %In that case, 
In order to do so, one just has to use the vectorial setpoint quantity $\bm{y} :=  (F_{pP}, F_{wG})$ in~\eqref{eq:goal-fct-setpoint-tracking} and minimize the resulting goal function w.r.t.~the vectorial control quantity $(F_{fB}, T)$. As has been remarked above (Section~\ref{sect:pi-results}), the feedback-based approach, by contrast, has its limitations when applied to several control channels simultaneously. %applying several proportional-integral controllers at the same time (simultaneously), by contrast, has its limitations

%\begin{figure}[htbp]
%\centering
%\includegraphics[width=\textwidth]{setpoint-tracking-Ft_G-by-opt-control.jpg}
%\caption{Setpoint tracking for $F_{wG}$ by optimal control}
%\label{fig:setpt-Fw_G-opt-contr}
%\end{figure}

\section{Summary and discussion}
%\section{Conclusion and perspectives}

In this paper, we dealt with the simulation and optimization of a prototypical process from chemical engenieering, the Williams-Otto process, using Pyomo, a powerful and versatile %and user-friendly %high-level 
dynamic simulation and optimization %environment
software package based on Python. 
We discussed two typical optimal control problems for the Williams-Otto process, namely waste minimization and yield maximization. 
In order to set up and solve such optimal control problems in Pyomo, one basically only has to complement the dynamic process equations by the appropriate goal function and constraints. 
As one would expect, waste minimization and yield maximization are conflicting goals and we briefly discussed the combined optimization of these conflicting goals by a weighted-sum approach. In view of the non-convexity of the problem, it would be natural to further investigate this multicriteria optimization problem, however  (Pareto front approximation and navigation ~\cite{Te14}, \cite{BoBu+14}, \cite{BuBo+14}, for instance). A first study in that direction, for the reactor subsystem of the Williams-Otto process, is undertaken in~\cite{LoVaHoDi12}. %A first mulitcriteria-optimization study for the reactor subsystem of the Williams-Otto process is undertaken in~\cite{LoVaHoDi12}. %We leave that to future research. 
%\smallskip
%
Aside from process optimization, we also discussed two approaches to the setpoint tracking problem:
\begin{itemize}
\item one standard approach which is based on feedback control, more precisely, on proportional-integral control, and
\item one approach that is based on optimal control.
\end{itemize}
We saw that with the optimization-based approach, setpoints are typically reached more quickly and more reliably. 
In order to perform setpoint tracking in Pyomo by means of proportional-integral control, one basically only has to complement the dynamic process equations by the appropriate feedback relations. 
It would be natural and interesting to further  investigate setpoint tracking using yet another approach, namely model-predictive control~\cite{Al}, \cite{GrPa}, \cite{FiAlBi07}, \cite{AnDi+19} which combines the feedback- and optimization-based approaches from above. 
Just like the two approaches disscussed above, a model-predictive control scheme for the Williams-Otto process could be implemented in a quite  straightforward manner with Pyomo. (See \cite{Am11} for a discussion of  model-predictive control for the reactor subsystem of the Williams-Otto process.) 
%As far as we know, model-predictive control of the Williams-Otto process has not been thoroughly studied so far (although there are studies~\cite{Am11} discussing model-predictive control of the reactor subsystem of the Williams-Otto process). %We leave that to future research.    

%Summarizing, %Concluding/And finally, 
%we saw that simulation, setpoint tracking, and optimization of dynamic processes like the Williams-Otto process can be conveniently done using the Pyomo environment. 
%And finally, throughout our discussions, we saw the strength, versatility, and easy usability of the software environment Pyomo in dealing with simulation, setpoint tracking, and optimization of dynamic processes like the Williams-Otto process. 

\section{Conclusion}

In this paper, we illustrated the advantages that a high-level software package (Pyomo) can provide for the chemical engineer. We paradigmatically solved a number of increasingly complex dynamic simulation and optimization problems based on the Williams-Otto process. We performed waste and yield optimization. Additionally, we contrasted the use of proportional-integral  controllers with an open-loop optimization approach for setpoint tracking.

\section*{Acknowledgment}

We thank Paul Beckermann for his assistance in preparing Section~\ref{sect:pyomo}. %contribution to Section~\ref{sect:pyomo}. 

\begin{small}

\end{small}


\begin{thebibliography}{}
%\bibliographystyle{alpha}
\bibitem{Al} F. Allg\"ower, A. Zheng (editors): Nonlinear model-predictive control. Birkh\"auser (2000)

\bibitem{Am11} R. Amrit: Optimizing process economics in model-predictive control. Ph.D. thesis, University of Wisconsin, Madison (2011)

\bibitem{AnDi+19} J.A.E. Andersson, J. Gillis, G. Horn, J.B. Rawlings, M. Diehl: \emph{CasADi: a software framework for nonlinear optimization and optimal control.} Math. Prog. Comp. \textbf{11} (2019), 1–36

\bibitem{AsMa14} F. Assassa, W. Marquardt: \emph{Dynamic optimization using adaptive direct multiple shooting.} Comp. Chem. Eng.~\textbf{60} (2014), 242-259 

\bibitem{BhBi96} T. Bhatia, L. T. Biegler: \emph{Dynamic optimization in the design and scheduling of multiproduct batch plants.} Ind. Eng. Chem. Res. \textbf{35} (1996), 2234

\bibitem{Biegler10} L.T. Biegler: Nonlinear programming. Concepts, algorithms, and applications to chemical processes. MOS-SIAM Series on Optimization (2010)

\bibitem{Biegler97}  L.T. Biegler, I.E. Grossmann, A.W. Westerberg: Systematic methods of chemical process design. Prentice Hall (1997)

\bibitem{BoPl84} H.G. Bock, K.J. Plitt: \emph{A multiple shooting algorithm for direct solution of optimal control problems.} IFAC Proceedings Volumes \textbf{17} (1984), 1603-1608 

\bibitem{BoBu+14} M. Bortz, J. Burger, N. Asprion, S. Blagov, R. Böttcher, U. Nowak, A. Scheithauer, R. Welke, K.-H. Küfer, H. Hasse: \emph{Multi-criteria optimization in chemical process design and decision support by navigation on Pareto sets.} Comp. Chem. Eng. \textbf{60} (2014), 354-363

\bibitem{BoHe+19} M. Bortz, R. Heese, A. Scherrer, T. Gerlach, T. Runowski: \emph{Estimating mixture properties from batch
distillation using semi-rigorous and rigorous
models.} In A.A. Kiss, E. Zondervan, R. Lakerveld, L. \"Ozkan (editors): 
Proceedings of the 29th European Symposium on Computer Aided Process Engineering. Elsevier (2019)

\bibitem{BressanPiccoli} A. Bressan, B. Piccoli: Introduction to the mathematical theory of control. AIMS Series on Applied Mathematics (2007)

\bibitem{BuBo+14} J. Burger, N. Asprion, S. Blagov, R. Böttcher, U. Nowak, M. Bortz, R. Welke, K.-H. Küfer, H. Hasse: \emph{Multi‐objective optimization and decision support in process engineering – Implementation and application.} Chem. Ing. Techn. \textbf{86} (2014), 1065-1072

\bibitem{CaKw16} A. Oliveira Cardoso, W.H. Kwong: \emph{Williams-Otto plant control based on production planning associated to coordinated decentralized optimization and plantwide control techniques.} J. Chem. Chem. Eng. \textbf{10} (2016), 77-89

\bibitem{CeBi98} A. Cervantes, L.T. Biegler: \emph{Large-scale DAE optimization using simultaneous nonlinear programming formulations.} AIChE J. \textbf{44} (1998), 1038

\bibitem{DiBo+06} M. Diehl, H.G. Bock, H. Diedam, P.-B. Wieber: Fast direct multiple shooting algorithms for optimal robot control. In: Fast motions in biomechanics and robotics, Lecture Notes in Control and Information Sciences \textbf{340} (2006), 65-93

\bibitem{FiAlBi07} R. Findeisen, F. Allg\"ower, L.T. Biegler (editors): Assessment and future directions of nonlinear model predictive control. Springer  (2007)


%\bibitem{GmehlingOnken} M. Baerns, A. Behr, A. Brehm, J. Gmehling, H. Hofmann, U. Onken,  A. Renken, K.-O. Hinrichsen, R. Palkovits: Technische Chemie. 2nd edition, Wiley (2013)

\bibitem{GrKrRa01} M. Gr\"otschel, S. Krumke, J. Rambau (editors): Online optimization of large systems. Springer (2001)

\bibitem{GrPa} L. Gr\"une, J. Pannek: Nonlinear model predictive control. Theory and algorithms. 2nd edition, Springer (2017)

\bibitem{HaMa12} R. Hannemann-Tamas, W. Marquardt: \emph{How to verify optimal controls computed by direct shooting methods? -- A tutorial.} J. Proc. Contr. \textbf{22} (2012), 494-507

\bibitem{Ha+11} W.E. Hart, J.-P. Watson, D.L. Woodruff: \emph{Pyomo:  modeling and solving mathematical programs in Python.} Math. Program. Comp. \textbf{3} (2011), 219-260 

\bibitem{Pyomo} W.E. Hart, C.D. Laird, J.-P. Watson, D.L. Woodruff, G.A. Hackebeil, B.L. Nicholson, J.D. Siirola: Pyomo -- Optimization modeling in Python. 2nd edition, Springer (2017) 

\bibitem{HoFeDi11} B. Houska, H.J. Ferreau, M. Diehl: \emph{ACADO toolkit -- An open-source framework for automatic control and dynamic optimization.} 
Optim. Control Appl. Meth. \textbf{32} (2011), 298–312 


%\bibitem{Levenspiel}  O. Levenspiel: Chemical reaction engineering. 3rd edition, Wiley (1999)

\bibitem{LoVaHoDi12} F. Logist, M. Vallerio, B. Houska, M. Diehl, J. Van Impe: \emph{Multi-objective optimal control of chemical processes using ACADO toolkit.} Comp. Chem. Eng.~\textbf{37} (2012), 191-199

\bibitem{Ni+18} B.L. Nicholson, J.D. Siirola, J.-P. Watson, V.M. Zavala, L.T. Biegler: \emph{pyomo.dae: a modeling and automatic discretization framework for optimization with differential and algebraic equations.} Math. Progr. Comp. \textbf{10} (2018), 187–223 

\bibitem{O'Dwyer12} A. O'Dwyer: \emph{An overview of tuning rules for the PI and PID continuous-time control of time-delayed single-input, single-output (SISO) processes.} In R. Vilanova, A. Visioli: PID control in the third millennium, Springer (2012)

\bibitem{Sk06} S. Skogestad: \emph{Tuning for smooth PID control with acceptable disturbance rejection}. Ind. Eng.
Chem. Res. \textbf{45} (2006), 7817–7822 

\bibitem{SkGr12} S. Skogestad, C. Grimholt: \emph{The SIMC method for smooth PID controller tuning.} In R. Vilanova, A. Visioli: PID control in the third millennium, Springer (2012)

\bibitem{Sontag} E.D. Sontag: Mathematical control theory. Deterministic finite-dimensional systems. 2nd edition, Springer (1998)

\bibitem{Te14} K. Teichert: A hyperboxing Pareto approximation method applied to radiofrequency ablation treatment planning. Ph.D. thesis, Technical University of Kaiserslautern (2014)

\bibitem{TjBi91} I.-B. Tjoa, L.T. Biegler: \emph{Simultaneous solution and optimization strategies for parameter estimation of differential-algebraic equation systems.} Ind. Eng. Chem.
Res. \textbf{30} (1991), 376–385

\bibitem{ViVi12} R. Vilanova, A. Visioli: PID control in the third millennium. Springer (2012)

\bibitem{WaBi05a} A. W\"achter, L.T. Biegler: \emph{Line Search Filter Methods for Nonlinear Programming: Motivation and Global Convergence.} SIAM J. Optim. \textbf{16} (2005), 1–31

\bibitem{WaBi05b} A. W\"achter, L.T. Biegler: \emph{Line Search Filter Methods for Nonlinear Programming: Local Convergence.} SIAM J. Optim. \textbf{16} (2005), 32-48

\bibitem{WaBi06} A. W\"achter, L.T. Biegler: \emph{On the implementation of an interior-point filter line-search algorithm for large-scale nonlinear programming.} Math. Progr. \textbf{106} (2006), 25–57

\bibitem{WiOt60} T. Williams, R. Otto: \emph{A generalized chemical processing model for the investigation of computer control.} AIEE Trans. \textbf{79} (1960), 458-473

\end{thebibliography}
\end{document}